% Paper: 	" On holomorphic domination, I "
%		
% Date:   	November 27, 2008.
% Author: 	(Mr.) Imre PATYI
%         	Department of Mathematics and Statistics
%         	Georgia State University
%		Atlanta, GA 30303-3083
%         	USA
% Email:        imrepatyi@speedpost.net
% Tel. 404-413-6450
% Fax. 404-413-6403
% Remarks on the file:

% 0.  This is an AmS-TeX file, not AmS-LaTeX.
% 1.  There are many abbreviations and a few macros.
% 2.  If any difficulty arises with the file or the
%    content of this manuscript, please email the author
%    at "ipatyi@gsu.edu".
% 3.  This file contains 2408 lines.

% PAPER BEGIN
% ~~~~~~~~~~~~~~~~~~~~~~~~~~~~~~~~~~~~~~~~~~~~~~~~~~~~~~~~~~~~~~~~~~~~~~~~~~~~~
% ~~~~~~~~~~~~~~~~~~~~~~~~~~~~~~~~~~~~~~~~~~~~~~~~~~~~~~~~~~~~~~~~~~~~~~~~~~~~~
% ~~~~~~~~~~~~~~~~~~~~~~~~~~~~~~~~~~~~~~~~~~~~~~~~~~~~~~~~~~~~~~~~~~~~~~~~~~~~~
% ~~~~~~~~~~~~~~~~~~~~~~~~~~~~~~~~~~~~~~~~~~~~~~~~~~~~~~~~~~~~~~~~~~~~~~~~~~~~~

%\input /math/faculty3/ipatyi/tex/amstex.tex

\input amstex.tex

\magnification=\magstep1
\hsize=5.5truein
\vsize=9truein
\hoffset=0.5truein
\parindent=10pt
\newdimen\nagykoz
\newdimen\kiskoz
\nagykoz=7pt
\kiskoz=2pt
\parskip=\nagykoz
\baselineskip=12.7pt

%\raggedbottom

\loadeufm \loadmsam \loadmsbm

\font\vastag=cmssbx10
\font\drot=cmssdc10
\font\vekony=cmss10
\font\vekonydolt=cmssi10
\font\cimbetu=cmssbx10 scaled \magstep1
\font\szerzobetu=cmss10

\font\scVIII=cmcsc8
\font\rmVIII=cmr8
\font\itVIII=cmti8
\font\bfVIII=cmbx8
\font\ttVIII=cmtt8

\def\cim#1{{\centerline{\cimbetu#1}}}
\def\szerzo#1{{\vskip0.3truein\centerline{\szerzobetu#1}}}
\def\alcim#1{{\medskip\centerline{\vastag#1}}}
\def\tetel#1#2{{{\drot#1}{\it\szukebb~#2\tagabb}}}
\long\def\biz#1#2{{{\vekony#1} #2}}
\def\kiemel#1{{\vekonydolt#1\/}}
\long\def\absztrakt#1#2{{\vskip0.4truein{\vekony#1} #2\vskip0.5truein}}
\def\szukebb{\parskip=\kiskoz}
\def\tagabb{\parskip=\nagykoz}
\def\vonal{{\vrule height 0.2pt depth 0.2pt width 0.5truein}}

\def\CC{{\Bbb C}}

\def\bmfd{{Banach manifold}}
\def\bsmfd{{Banach submanifold}}

\def\bvbdl{{Banach vector bundle}}

\def\re{\text{Re}}

\def\cln{{:}\,}

\def\Oa{{\Omega}}
\def\oa{{\omega}}

\def\cts{{continuous}}
\def\bdd{{bounded}}

\def\idml{{infinite dimensional}}

\def\fdml{{finite dimensional}}

\def\pscx{{pseudoconvex}}

\def\bd{{\partial}}

\def\bap{{bounded approximation property}}
\def\sbs{{Schauder basis}}
\def\ubs{{unconditional basis}}

\def\psh{{plurisubharmonic}}
\def\pshdom{{\psh\ domination}}
\def\holodom{{\holo\ domination}}
\def\domentfb{{dominated by entire \fns\ with values in \bspc{}s}}

\def\st{{such that}}

\def\subsub{\subset\!\subset}
\def\delbar{{\bold{\bar\partial}}}

\def\<{{\langle}}
\def\>{{\rangle}}

\def\RR {{\Bbb R}}

\def\OO {{\Cal O}}

\def\ro{\varrho}

\def\epsz{\varepsilon}
\def\fii{\varphi}
\def\fn{func\-tion}
\def\fns{func\-tions}
\def\holo{hol\-o\-mor\-phic}

\def\cpx{complex}

\def\nbd{neighbor\-hood}

\def\bspc{Banach space}

\def\la{\lambda}

\def\aa{{\alpha}}
\def\ba{{\beta}}

\def\ta{{\theta}}

\def\WLOG{Without loss of generality}

\def\za{\zeta}

\def\Prop{Proposition}
\def\p#1{{\Prop~#1}}

\def\Th{{Theorem}}
\def\th{{theorem}}

\def\t#1{{\Th~#1}}

% today.tex: Macro to print today's date
%--
% Author: John W. Shipman, NM Tech Computer Center,
%   Socorro, NM 87801; john at nmt.edu
%--
% EXPORTED FUNCTIONS:
%   \today: Outputs today's date as ``yyyy-mm-dd''
%   \now: Outputs the current time as ``hh:dd''
%   \timestamp: \today, plus one space, plus \now
%--
\newcount\minute    % Current minute within the hour
\newcount\hour      % Current hour (24-hour type)
\newcount\hourMins  % Temporary for taking hour modulo 60
%
% - - -   \ n o w   - - -
%
\def\now%
{% Displays today's time as ``hh:mm''
%  The \time macro gives the minutes since midnight.  Compute
%  the whole hours by dividing this by 60, then find the
%  minute by effectively taking the minutes modulo 60.
%
  \minute=\time    % Number of minutes since midnight
  \hour=\time \divide \hour by 60 % Get hours
  \hourMins=\hour \multiply\hourMins by 60
  \advance\minute by -\hourMins % Hours modulo 60
  \zeroPadTwo{\the\hour}:\zeroPadTwo{\the\minute}%
}% --- \now ---
%
% - - -   \ t i m e s t a m p   - - -
%
\def\timestamp%
{% Displays ``yyyy-mm-dd hh:mm''
  \today\ \now
}% --- \timestamp ---
%
% - - -   \ t o d a y   - - -
%
\def\today%
{% Displays today's date and time as ``yyyy-mm-dd hh:mm''
  \the\year-\zeroPadTwo{\the\month}-\zeroPadTwo{\the\day}%
}% --- \today ---
%
% - - -   z e r o P a d T w o   - - -
%
\def\zeroPadTwo#1%
{% Left zero pad of the argument to 2 digits.  The argument
%  should be a number between 1 and 99.  This macro outputs
%  a `0' if the argument is less than ten, then it outputs
%  the argument.
%
  \ifnum #1<10 0\fi    % Conditionally output a zero
  #1%    Then output the argument
}% --- \zeroPadTwo ---

%{\ttVIII\jobname, \timestamp}
{\phantom.}
\vskip0.5truein
\cim{ON HOLOMORPHIC DOMINATION, I}
%\vskip0.2truein
%\cim{CERTAIN COMPLEX BANACH MANIFOLDS}
\szerzo{Imre Patyi\plainfootnote{$\,{}^1$}{\rmVIII 
	Supported in part by an NSF Grant DMS-0600059.}}
\absztrakt{ABSTRACT.}{Let $X$ be a separable Banach space and $u{:}\,X\to\RR$
	locally upper bounded.
	 We show that there are a Banach space $Z$ and a holomorphic function $h{:}\,X\to Z$
	with $u(x)<\|h(x)\|$ for $x\in X$.
	 As a consequence we find that the sheaf cohomology group $H^q(X,\OO)$
	vanishes if $X$ has the bounded approximation property (i.e., $X$ is a
	direct summand of a Banach space with a Schauder basis), 
	$\OO$ is the sheaf of germs of
	holomorphic functions on $X$, and $q\ge1$.
	 As another consequence we prove that if $f$ is a $C^1$-smooth 
	$\delbar$-closed $(0,1)$-form on the space $X=L_1[0,1]$ of summable functions,
	then there is a $C^1$-smooth function $u$ on $X$
	with $\delbar u=f$ on $X$.

	

%Running head: On holomorphic domination, I, I. Patyi
\smallskip
{\it\hfill Kedves Zoli \"ocs\'emnek, sz\"ulet\'esnapj\'ara.\plainfootnote{$\,{}^2$}{\rmVIII To my younger brother.}}
}

% Section numbering
\def\sA{{1}}
\def\sB{{2}}
\def\sC{{3}}
\def\sD{{4}}
\def\sE{{5}}
\def\sF{{6}}

% Proclaim numbering (theorems) \tAA etc.
% Display numbering (equations) \eAA etc.

\def\tAA{{\sA.1}}
\def\tAB{{\sA.2}}
\def\tAC{{\sA.3}}
\def\tAD{{\sA.4}}
\def\tAE{{\sA.5}}

\def\tCA{{\sC.1}}
\def\tCB{{\sC.2}}
\def\tCC{{\sC.3}}

\def\eCA{{\sC.1}}

\def\tDA{{\sD.1}}
\def\tDB{{\sD.2}}

\def\tEA{{\sE.1}}
\def\tEB{{\sE.2}}
\def\tEC{{\sE.3}}
\def\tED{{\sE.4}}
\def\tEE{{\sE.5}}

\def\tFA{{\sF.1}}
\def\tFB{{\sF.2}}

% Reference numbering
\def\rDZ{DZ}
\def\rDPV{DPV}
\def\rLA{L1}
\def\rLB{L2}
\def\rLC{L3}
\def\rLP{LP}
\def\rM{M}
\def\rP{P}
\def\rPA{Pt1}
\def\rPB{Pt2}

\alcim{\sA. INTRODUCTION.}

	 The ideas of \pshdom\ and \holodom\ along with some of their applications
	appeared in [\rLC] by Lempert.
	 Following him we say that \kiemel{\pshdom} is possible on a \cpx\ \bmfd\ $M$ 
	if for every locally upper bounded $u\cln M\to\RR$ there is a \cts\ \psh\
	$\psi\cln M\to\RR$ with $u(x)<\psi(x)$ for all $x\in M$.
	 If $\psi$ can be taken in the form $\psi(x)=\|h(x)\|$ for $x\in M$, where
	$h\cln X\to Z$ is a \holo\ \fn\ to a \bspc\ $Z$, then we say that 
	\kiemel{\holodom} is possible in $M$.

	 One tool to achieve \holodom\ is the following Runge approximation property
	of a \bspc\ $X$.

\tetel{Hypothesis~\tAA.}{{\rm[\rLC, Hypothesis~1.5]}
	 There is a constant\/ $0<\mu<1$ \st\ if $Z$ is any \bspc, $\epsz>0$, and
	$f\cln B_X\to Z$ is \holo\ on the open unit ball $B_X$ of $X$, then
	there is a \holo\ \fn\ $g\cln X\to Z$ with\/ $\|f(x)-g(x)\|<\epsz$ for\/ $\|x\|<\mu$.
}

	 Lempert and Meylan proved the following \th\ involving the above.

\tetel{\t\tAB.}{{\rm(a) (Lempert, [\rLB])} If $X$ is a \bspc\ with an \ubs, then
	Hypothesis~\tAA\ above holds for $X$.
\vskip0pt
	{\rm(b) (Meylan, [\rM])} If $X$ is a \bspc\ with an unconditional \fdml\
	Schauder decomposition, then Hypothesis~\tAA\ holds for $X$.
\vskip0pt
	{\rm(c) (Lempert, [\rLC])} 
	If $X$ is a \bspc\ with a \sbs\ (or a direct summand of one)
	and Hypothesis~1.1 holds for $X$, then \holodom\ is possible in every
	\pscx\ open subset of $X$.
}

	 Our main goal in this paper is to find a route to \holodom\ that bypasses
	Hypothesis~\tAA\ above.
	 Our main results are \Th{}s \tAC, \tAD, \tAE, and \tFA\ below.

\tetel{\t\tAC.}{If $X$ is a separable \bspc, then \holodom\ is possible\/
	{\rm(a)} in $X$, and\/ {\rm(b)} in every convex open\/ $\Oa\subset X$.
}

	 As a consequence of \t\tAC\ we get cohomology vanishing as follows.

\tetel{\t\tAD.}{Let $X$ be a \bspc\ with the \bap, $\Oa\subset X$ \pscx\ open,
	$M\subset\Oa$ a closed split \cpx\ \bsmfd\ of\/ $\Oa$, $S\to M$ a cohesive
	sheaf, $E\to\Oa$ a \holo\ \bvbdl, and $I\to\Oa$ the sheaf of germs of \holo\
	sections of $E$ over\/ $\Oa$ that vanish on $M$.
	 If \pshdom\ is possible in\/ $\Oa$ (which is guaranteed by \t\tAC\ if\/
	$\Oa\subset X$ is convex open), then the following hold.
\vskip0pt
	{\rm(a)} The cohesive sheaf $S\to M$ admits a complete resolution over $M$.
\vskip0pt
	{\rm(b)} The sheaf cohomology group $H^q(M,S)$ vanishes for all $q\ge1$.
\vskip0pt
	{\rm(c)} The sheaf $I$ is cohesive over\/ $\Oa$, $H^q(\Oa,I)=0$ for $q\ge1$,
	and any \holo\ section $f\in\OO(M,E)$ extends to a \holo\ section $F\in\OO(\Oa,E)$
	with $F(x)=f(x)$ for $x\in M$.
\vskip0pt
	{\rm(d)} If\/ $\Oa\subset X$ is convex open, then $E$ is \holo{}ally trivial
	over\/ $\Oa$.
}

	 As a consequence of \t\tAD\ we get the following \t\tAE\ on the $\delbar$-equation.

\tetel{\t\tAE.}{Let $X$ be an ${\Cal L}_1$-space with the \bap\ (e.g., $X=L_1[0,1]$),
	$\Oa\subset X$ \pscx\ open, $E\to\Oa$ a \holo\ \bvbdl, and $f\in C^1_{0,1}(\Oa,E)$
	a $C^1$-smooth $\delbar$-closed\/ $(0,1)$-form with values in $E$.
	 If \pshdom\ is possible in\/ $\Oa$ (which is guaranteed by \t\tAC\ if\/
	$\Oa\subset X$ is convex open), then there is a $C^1$-smooth section
	$u\in C^1(\Oa,E)$ of $E$ with $\delbar u=f$ over\/ $\Oa$.
}

	 Our strategy is to imitate the relevant parts of [\rLC] closely, but
	refrain from using Runge approximation for \fns\ un\bdd\ on balls.
	 The reader is assumed to have a copy of [\rLC] along side this paper.	
	 In our \S\S\,\sB-\sD\ we adopt without comment the notation of [\rLC, \S\S\,2-4].

\alcim{\sB. BACKGROUND.}

	 In this section we recall some material useful later.
	 The paper [\rLC] uses a particular exhaustion $\Oa_N\<\aa\>$, $N\ge1$,
	of any \pscx\ open subset $\Oa$ of any \bspc\ $X$ with a bimonotone \sbs,
	and there are numerous other sets used there to help out with the analysis
	of the said exhaustion.
	 In our case all the sets involved will be convex open in $X$ or in the
	span of finitely many of its basis vectors.
	 The \idml\ ones among the sets that we need are all of the form $D\times B$,
	where $D$ is a convex open set in the span of the first few basis vectors
	and $B$ is a ball in the closed span of the rest of the basis vectors.
	 As we shall need very little of the properties of the many sets discussed
	in [\rLC] we just help ourselves directly to the results there and
	skip any of their details (even their definitions) here.

	 In a \bspc\ $X$, put $B_X(x_0,r)=\{x\in X\cln\|x-x_0\|<r\}$ for the open ball
	of radius $r$ centered at $x_0\in X$, and write $B_X=B_X(0,1)$ for the unit ball.
	 Denote by $\OO(M_1,M_2)$ the set of \holo\ \fns\ $M_1\to M_2$ from one \cpx\
	\bmfd\ $M_1$ to another $M_2$.

	 Let $X$ be a \bspc, $A\subset X$, and $u\cln A\to\RR$.
	 We say that $u$ can be \kiemel{\domentfb} on $A$ if there are a \bspc\ $Z$
	and an entire \holo\ \fn\ $h\in\OO(X,Z)$ with $u(x)\le\|h(x)\|$ for all $x\in A$.

	 If $T$ is any set, then denote by $\ell_\infty(T)$
	the \bspc\ of \bdd\ \fns\ $f\cln T\to\CC$ with the sup norm 
	$\|f\|=\sup\{\|f(t)\|: t\in T\}$.

\alcim{\sC. DOMINATION ON THE WHOLE SPACE.}

	 In this section we show that if a \fn\ $u$ can be dominated on every ball
	of a fixed radius, then $u$ can be dominated on the whole space as well.

	 Let $X$ be \bspc\ with a \sbs.
	 Fix the norm and the \sbs\ of $X$ so as to make a bimonotone \sbs\ of $X$.
	 Fix $N\ge1$ and write $\pi$ for the Schauder projection onto the span
	of the first $N+1$ basis vectors, $\ro=1-\pi$ for the complementary projection,
	and $Y=\ro X$ for the complementary space.

\tetel{\p\tCA.}{If $X$ is a \bspc\ with a bimonotone \sbs, $0<R<\infty$, 
	$u\cln X\to[1,\infty)$ is \cts, and $u$ can be \domentfb\ on every ball
	$B_X(x_0,R)$ of radius $R$ and centered at any $x_0\in X$, then
	u can be \domentfb\ on $X$.
}

	 The proof of \p\tCA\ will occupy us for a while.

\tetel{\p\tCB.}{{\rm(Cf.\ [\rLC, Lemma 4.1])}
	 Let $A_2\subsub A_3$ be relatively open \bdd\ convex subsets of 
	$\pi(X)\cong\CC^{N+1}$, $A_1$ a compact convex subset of $A_2$, and
	$0<r_1<r_2<r_3<\infty$ constants.
	 If $Z$ is a \bspc\ and $g\in\OO(X,Z)$ is an entire \fn, then there are
	a \bspc\ $W$ and an entire \fn\ $h\in\OO(X,W)$ with
\vskip0pt
	{\rm(i)} $\|h(x)\|_W\le1$ for $x\in A_1[r_1]$ and
\vskip0pt
	{\rm(ii)} $\|h(x)\|_W\ge\|g(x)\|_Z$ for $x\in A_3(r_3)\setminus A_2(r_2)$.
}

\biz{Proof.}{Consider the \bdd\ convex sets $H_1,H_2,H_3$ in $\pi(X)\times\CC\cong\CC^{N+2}$
	given by $H_1=\{(s,\la)\in A_1\times\CC:|\la|\le r_1\}$,
	$H_i=\{(s,\la)\in A_i\times\CC:|\la|<r_i\}$ for $i=2,3$.
	 Since $H_1$ is compact convex in $\CC^{N+2}$ there are a finite set $J$ and
	polynomials $\fii_j\in\OO(\pi(X)\times\CC)$ for $j\in J$ \st\
	$|\fii_j(s,\la)|\le\frac14$ for $(s,\la)\in H_1$ and for every 
	$(s,\la)\in H_3\setminus H_2$ there is a $j\in J$ with $|\fii_j(s,\la)|\ge4$.
	 Denote by $L=\overline{B_{Y^*}}$ the set of all linear functionals $l\in Y^*$
	with $\|l\|\le1$, and by $V=\ell_\infty(L\times J)$.
	 Define $\fii\in\OO(X,V)$ by $\fii(x)(l,j)=\fii_j(\pi x,l\ro x)$ for $x\in X$,
	$l\in L$, and $j\in J$.

	 The rest of the proof of \p\tCB\ is the same word for word as that of
	[\rLC, Lemma 4.1] starting with ``Going back'' near [\rLC, (4.1)].
}

\tetel{\p\tCC.}{{\rm(Cf.\ [\rLC, \p4.2])} Let\/ $0<\mu<1$, $N\ge1$,
	and $2^4\ba<\aa<2^{-8}\mu$.
	 If $Z$ is a \bspc\ and $g\in\OO(X,Z)$ is an entire \fn, then there are a
	\bspc\ $W$ and an entire \fn\ $h\in\OO(X,W)$ \st\
\vskip0pt
	{\rm(i)} $\|h(x)\|_W\le1$ for $x\in\Oa_N\<\ba\>$ and
\vskip0pt
	{\rm(ii)} $\|h(x)\|_W\ge\|g(x)\|_Z$ for $x\in\Oa_{N+1}\<\aa\>\setminus\Oa_N\<\aa\>$.
}

\biz{Proof.}{
	 In \p\tCC\ the sets $\Oa_N\<\ba\>$, etc, refer to those constructed in
	[\rLC, \S3] for $\Oa=X$.
	 \p\tCC\ follows from \p\tCB\ in the same way as [\rLC, \p4.2]
	does from [\rLC, Lemma 4.1] only more simply.
}

\biz{Proof of \p\tCA.}{On replacing $u$ by $u(Rx/2)$ we may assume that $R=2$.
	 Let $\Oa=X$, fix $0<\mu<1$ and $0<\aa<2^{-8}\mu$.
	 First, we construct a \bspc\ $Z_N$ and an entire \fn\ $g_N\in\OO(X,Z_N)$
	for each $N\ge1$.
	 The set $A=\overline{\Oa_N\<\aa\>}\cap\pi_N(X)$ is compact and if $t\in A$,
	then $\Oa_N\<\aa\>\cap\pi_N^{-1}(t)\subset B_X(t,\aa)$.
	 Hence $t$ has an open \nbd\ $U\subset\pi_N(X)$ with
	$\Oa_N\<\aa\>\cap\pi_N^{-1}(U)\subset B_X(t,2\aa)$.
	 Therefore
$$
	\Oa_N\<\aa\>\subset\bigcup_{t\in T} B_X(t,2\aa)
\tag\eCA
 $$
 	for some finite $T\subset A$.
	 Let $B_t=B_X(t,2\aa/\mu)$, the radius of which is less than $2$.
	 By our assumption that $u$ can be \domentfb\ on $B_X(x_0,2)$ for every $x_0\in X$,
	there are a \bspc\ $V_t$ and an entire \fn\ $f_t\in\OO(X,V_t)$ with
	$u(x)\le\|f_t(x)\|_{V_t}$ for $x\in B_t$, $t\in T$.
	 Let $Z_N$ be the $\ell_\infty$-sum of the finitely many \bspc{}s $V_t$ for $t\in T$
	and $g_N\in\OO(X,Z_N)$ the map whose components are the $f_t$ for $t\in T$.
	 We see from $(\eCA)$ that $u(x)\le\|g_N(x)\|_{Z_N}$ for $x\in\Oa_N\<\aa\>$.

	 The rest of the proof of \p\tCA\ is the same as that of [\rLC, \p2.1]
	starting with ``In the second step'' on page 368 there.
}

\alcim{\sD. DOMINATION ON A BALL.}

	 In this section we show that if a \fn\ $u$ can be dominated on every ball
	of half the radius of a ball $B$ and centered at any point of $B$, then
	$u$ can be dominated on $B$ itself.

\tetel{\p\tDA.}{If $X$ is a \bspc\ with a bimonotone \sbs, $0<R<\infty$,
	$u\cln X\to[1,\infty)$ is \cts, and $u$ can be \domentfb\ on every ball
	$B_X(x_0,R/2)$ of radius $R/2$ and centered at any $x_0\in B=B_X(y_0,R)$,
	then there is \cts\ \fn\ $\tilde u\cln X\to[1,\infty)$ \st\
	$\tilde u(x)\le u(x)$ for all $x\in X$, $\tilde u(x)=u(x)$ for $x\in B$,
	and $\tilde u$ can be \domentfb\ on every ball $B_X(x_0,R/8)$ of radius $R/8$
	centered at any $x_0\in X$.
}

\biz{Proof.}{Let $\chi\cln[0,\infty)\to[0,1]$ be a cutoff \fn
$$
	\chi(t)=\cases
		1&\text{\phantom{ifpea}}0\le t\le R\cr
		1-\frac4R(t-R)&\text{if\phantom{pea}}R\le t\le\frac54R\cr
		0&\text{\phantom{ifpea}}t\ge\frac54R
		\endcases,
 $$
 	and define $\tilde u$ by 
	$\tilde u(x)=\chi(\|x-y_0\|)u(x)+1-\chi(\|x-y_0\|)$ for $x\in X$.
	 As $\tilde u(x)-u(x)=(1-\chi(\|x-y_0\|))(1-u(x))\le0$, being the product
	of a nonnegative number by a nonpositive number, we get that $\tilde u(x)\le u(x)$
	for all $x\in X$.
	 Hence $\tilde u$ can be \domentfb\ on any set on which $u$ can.

	 If $x_0\in X$ satisfies that $\|x_0-y_0\|\ge\frac{11}{8}R$, then
	$B_X(x_0,\frac18R)$ lies outside $B_X(y_0,\frac54R)$ since
	the distance $\|x_0-y_0\|$ of their centers exceeds the sum of their radii
	$\frac54R+\frac18R=\frac{11}8R$.
	 Hence $\tilde u=1$ on $B_X(x_0,\frac18R)$, and so $\tilde u$ can be \domentfb\
	on $B_X(x_0,\frac18R)$.

	 If $\|x_0-y_0\|<R$, then $x_0\in B_X(y_0,R)$ and
	$B_X(x_0,\frac18R)\subset B_X(x_0,\frac12R)$.

	 If $R\le\|x_0-y_0\|<\frac{11}{8}R$, then choose a value $0<R'<R$ with
	$\|x_0-y_0\|<\frac{11}{8}R'$, and let $z_0=y_0+R'\frac{x_0-y_0}{\|x_0-y_0\|}$.
	 Then $\|z_0-x_0\|=R'<R$ so $z_0\in B_X(y_0,R)$ and we claim that
	$B_X(x_0,\frac18R)\subset B_X(z_0,\frac12R)$.
	 To that end we must show that the distance $\|z_0-x_0\|$ of the centers is
	less than the difference of the radii, i.e., 
	$\|z_0-x_0\|<\frac12R-\frac18R=\frac38R$.
	 Indeed, $\|z_0-x_0\|=\|y_0-x_0+R'\frac{x_0-y_0}{\|x_0-y_0\|}\|=\|x_0-y_0\|-R'<
	\frac{11}{8}R'-R'=\frac38R'<\frac38R$.
	 The proof of \p\tDA\ is complete.
}

\tetel{\p\tDB.}{If $X$ is a \bspc\ with a bimonotone \sbs, $0<R<\infty$,
        $u\cln X\to[1,\infty)$ is \cts, and $u$ can be \domentfb\ on every ball
	$B_X(x_0,R/2)$ of radius $R/2$ centered at any $x_0\in B=B_X(y_0,R)$,
	then $u$ can be \domentfb\ on the ball $B$.
}

\biz{Proof.}{\p\tDA\ gives us a $\tilde u$ that can be \domentfb\ on every ball
	of radius $R/8$ in $X$.
	 \p\tCA\ gives us a \bspc\ $Z$ and an entire \fn\ $h\in\OO(X,Z)$ with
	$\tilde u(x)\le\|h(x)\|$ for all $x\in X$.
	 As $u(x)=\tilde u(x)\le\|h(x)\|$ for $x\in B$, the proof of \p\tDB\
	is complete.
}

\alcim{\sE. PREPARATION.}

	 This section is preparatory to the proofs of \Th{}s \tAC, \tAD, and \tAE.

	 Recall the following theorem of Pe\l{}czy\'nski's.

\tetel{\t\tEA.}{{\rm(Pe\l{}czy\'nski, [\rP])} A \bspc\ $X$ has the \bap\ if
	and only if $X$ is isomorphic to a direct summand of a \bspc\ $Y$ with
	a \sbs, i.e., there are a \bspc\ $Y$ with a \sbs\ and a direct decomposition
	$Y=Y_1\oplus Y_2$ of \bspc{}s \st\ $X\cong Y_1$.
}

	 In most of our proofs we can avoid dealing with \bspc{}s with the \bap,
	and only work with \bspc{}s with a \sbs.

\tetel{\p\tEB.}{Let $X$ be a \bspc\ with the \bap, and\/ $\Oa\subset X$ \pscx\ open.
	 If \pshdom\ is possible in\/ $\Oa$, then so is \holodom.
}

\biz{Proof.}{It is enough by \t\tEA\ to prove this when $X$ has a \sbs,
	in which case it follows from the argument of [\rLC], only more simply.
}

\tetel{\p\tEC.}{If $M_0$ is a closed \cpx\ \bsmfd\ of a \cpx\ \bmfd\ $M$, and \holodom\
	is possible on $M$, then \holodom\ is possible on $M_0$, too.
}

\biz{Proof.}{Let $u_0\cln M_0\to\RR$ be the locally upper bounded \fn\ to be dominated.
	 Define $u\cln M\to\RR$ by setting $u(x)=u_0(x)$ for $x\in M_0$ and $u(x)=0$
	otherwise.
	 Clearly, $u$ is locally upper bounded, $M_0$ being a closed subset of $M$.
	 If $Z$ is \bspc\ and $h\in\OO(M,Z)$ dominates $u$ on $M$, then the restriction
	$h_0$ of $h$ to $M_0$ is \holo\ and dominates $u_0$ in $M_0$.
	 The proof of \p\tEC\ is complete.
}

\tetel{\p\tED.}{If $M$ is a separable \cpx\ \bmfd\ that is bi\holo\ to a closed \bsmfd\
	of a \bspc\ $X$, then $M$ can be embedded in a separable \bspc\ as a closed \cpx\
	\bsmfd.
}

\biz{Proof.}{It is easy to see that the closed linear span of a separable subset of
	any \bspc\ is itself separable.
	 It is a standard \th\ that any separable \bspc\ is isomorphic to a closed
	linear subspace of the space $Y=C[0,1]$ of \cts\ \fns, and $Y$ has a \sbs.
	 Thus $M$ is bi\holo\ to a closed \cpx\ \bsmfd\ of $Y$,
	completing the proof of \p\tED.
}

\tetel{\p\tEE.}{Let $X$ be a \bspc, and\/ $\Oa\subset X$ open.
	 If one of\/ {\rm(a), (b), (c)} below holds, then\/ $\Oa$ is bi\holo\
	to a closed \cpx\ \bsmfd\  $M$ of a \bspc\ $Y$.
\vskip0pt
	{\rm(a)} $\Oa$ is convex.
\vskip0pt
	{\rm(b)} There is a direct decomposition $X=X_1\oplus X_2$ of \bspc{}s with
	$\dim_\CC(X_1)<\infty$, and $\Oa$ is of the form\/
	$\Oa=\{(x_1,x_2)\in D\times X_2\cln\|x_2\|<R(x_1)\}$, where $D\subset X_1$ is \pscx\
	(relatively) open, $R\cln D\to(0,\infty)$ is \cts\ and $-\log R$ is \psh\ on $D$.
\vskip0pt
	{\rm(c)} $\Oa$ is of the form\/ $\Oa=\{x\in\Oa'\cln\|f(x)\|<1\}$,
	where $\Oa'\subset X$ is open, the closure\/ $\overline\Oa\subset\Oa'$,
	and $f\in\OO(\Oa',Z_1)$ is \holo\ with values in a \bspc\ $Z_1$.
}

\biz{Proof.}{In each case we define a \bspc\ $Z$ and a \holo\ \fn\ $h\in\OO(\Oa,Z)$
	with $\liminf_{\Oa\ni x\to x_0}\|h(x)\|=\infty$ for each boundary point 
	$x_0\in\bd\Oa$.
	 Then the graph $M\subset Y=X\times Z$ of $h$ defined by 
	$M=\{(x,z)\in\Oa\times Z\cln z=h(x)\}$ does the job.

	 (a) (See also [\rPA, \p8.2].)
	 Assume as we may that $0\in\Oa$.
	 Let $p\cln X\to\RR$ be the Minkowski functional
	$p(x)=\inf\{\la>0\cln\frac{x}\la\in\Oa\}$ of the convex open set $\Oa$,
	and $K=\{\xi\in X^*\cln\re(\xi x)\le p(x)\text{\ for all\ }x\in X\}$.
	 Then $K\not=\emptyset$ is a convex subset of the dual space $X^*$ of $X$.
	 We endow $K$ with the weak star topology, in which $K$ is compact.

	 Let $Z=C([0,2\pi]\times K,\CC)$ be the usual \bspc\ with the sup norm,
	$g(t)=1/(1-t)$ for $t\in B_\CC$,
	and define for $x\in\Oa$ a \fn\ $h(x)\in Z$ by 
	$h(x)(\ta,\xi)=g(e^{i\ta}e^{\xi x-1})$.
	 Then $\|h(x)\|=\sup\{|h(x)(\ta,\xi)|\cln\ta\in[0,2\pi],\;\xi\in K\}
	\le\sup_{\ta,\xi} g(|e^{i\ta}e^{\xi x-1}|)\le 
	\sup_{\ta,\xi}g(e^{\re(\xi x)-1})\le g(e^{p(x)-1})$.
	 For every $x\in X$ the Hahn--Banach theorem gives a $\xi\in K$
	with $\re(\xi x)=p(x)$.
	 On choosing $\ta\in[0,2\pi]$ so that 
	$e^{i\ta}e^{\xi x-1}=e^{\re(\xi x)-1}=e^{p(x)-1}$, we find that
	$\|h(x)\|=g(e^{p(x)-1})$.
	 Hence, $h\in\OO(\Oa,Z)$, and $\|h(x)\|=1/(1-e^{p(x)-1})\to\infty$
	as $x\in\Oa$ tends to point $x_0\in X$ with $p(x_0)=1$, in particular,
	to any boundary point $x_0\in\bd\Oa$.

	 (b) Let $\oa=\{(x_1,\la)\in D\times\CC\cln x_1\in D,|\la|<R(x_1)\}$.
	 As $\oa$ is \pscx\ open in the \cpx\ Euclidean space $X_1\times\CC$,
	there is a proper \holo\ embedding $j\cln\oa\to\CC^N$ for $N$ high enough.
	 Let $K$ be the closed unit ball of the dual space $X_2^*$ of $X_2$
	endowed with the weak star topology, and for 
	$(x_1,x_2)\in\Oa$ define $h(x_1,x_2)\in Z=C(K,\CC)$ 
	(endowed with the sup norm) by
	$h(x_1,x_2)(\xi_2)=j(x_1,\xi_2x_2)$ for $\xi_2\in K$.
	 Note that $\|h(x_1,x_2)\|=\sup_{|\la|\le\|x_2\|}\|j(x_1,\la)\|\ge
	\|j(x_1,\|x_2\|)\|$ by the Hahn--Banach theorem,
	and the last tends to $\infty$ if $(x_1,\|x_2\|)$ tends to a boundary point
	of $\oa$, in particular, when $(x_1,x_2)$ tends in $\Oa$ to a boundary
	point of $\Oa$ in $X$.

	 (c) Let $K$ be the closed unit ball of the dual space $Z_1^*$ of $Z_1$
	endowed with the weak start topology, and $Z=C(K,\CC)$ with the sup norm.
	 For $x\in\Oa$ define $h(x)\in Z$ by $h(x)\za=g(\za f(x))$, where
	$\za\in K$ and $g(t)=1/(1-t)$ for $t\in B_\CC$ as in (a).
	 Then $\|h(x)\|=g(\|f(x)\|)$ for $x\in\Oa$ by the Hahn--Banach theorem,
	and $h\in\OO(\Oa,Z)$ is \holo.
	 If $x\in\Oa$ tends to a boundary point $x_0\in\bd\Oa$, then
	$x_0\in\overline{\Oa}\subset\Oa'$, hence $x_0\in\Oa'$ and $f(x)\to f(x_0)$,
	i.e., $\|f(x)\|\to\|f(x_0)\|=1$,
	and $\|h(x)\|=1/(1-\|f(x)\|)\to\infty$.

	 The proof of \p\tEE\ is complete.
}

\alcim{\sF. THE PROOFS OF THEOREMS \tAC, \tAD, AND \tAE.}

	 In this section we complete the proof of \Th{}s \tAC\ on \holodom,
	 \tAD\ on vanishing and \bvbdl{}s, and \tAE\ on the $\delbar$-equation.

\biz{Proof of \t\tAC(a).}{\WLOG\ we may assume by \t\tEA\ that $X$ has a bimonotone \sbs.
	 Let $u\cln X\to\RR$ be the locally upper bounded \fn\ to be dominated.
	 By paracompactness of $X$ there is a \cts\ \fn\ $u_1\cln X\to[1,\infty)$
	with $u(x)\le u_1(x)$ for $x\in X$.
	 Replacing $u$ by $u_1$, let us assume that $u\ge1$ is \cts\ on $X$.

	 Suppose for a contradiction that $u$ cannot be \domentfb\ on $X$.
	 The hypothesis of \p\tCA\ must then be false.	 
	 Hence there is a ball $B_0=B_X(x_0,1)$ on which $u$ cannot be \domentfb.
	 The hypothesis of \p\tDB\ must then also be false.
	 So there is a ball $B_1=B_X(x_1,1/2)$ with $x_1\in B_0$ \st\ $u$ cannot be
	\domentfb\ on $B_1$.
	 Again, the hypothesis of \p\tDB\ must be false and there is a ball
	$B_2=B_X(x_2,1/4)$ with $x_2\in B_1$ \st\ $u$ cannot be \domentfb\ on $B_2$.
	 Proceeding in this way we get a sequence of balls $B_n=B_X(x_n,1/2^n)$
	with $x_{n+1}\in B_n$ \st\ $u$ cannot be \domentfb\ on $B_n$ for $n\ge0$.

	 As $x_{n+1}\in B_n$ we see that $\|x_{n+1}-x_n\|<1/2^n$ and
	$\sum_{n=0}^\infty(x_{n+1}-x_n)$ is an absolutely convergent series in the \bspc\ $X$.
	 Thus there is a limit $x_n\to x\in X$ as $n\to\infty$.
	 Let $r>0$ be so small that $u$ is upper bounded on the ball $B_X(x,r)$.
	 Choose $n\ge0$ so large that $B_n\subset B_X(x,r)$.
	 Hence $u$ can be \domentfb\ on $B_n$ after all, being upper bounded there.
	 This contradiction completes the proof of (a).
}

\tetel{\t\tFA.}{{\rm(a)} If $M$ is as in \p\tED, then \holodom\ is possible in $M$.
\vskip0pt
	{\rm(b)} In particular, if $X$ is a separable \bspc, and\/ $\Oa\subset X$ open
	is as in \p\tEE, then \holodom\ is possible in\/ $\Oa$.
}

\biz{Proof.}{Part (a) follows from \t\tAC(a) via \p\tEC\ upon embedding $M$ 
	in $C[0,1]$ as a closed \cpx\ \bsmfd.
	 Part (b) follows from (a) by \p\tEE.
	 The proof of \t\tFA\ is complete, and as \t\tAC(b) is a special case
	of (b), the proof of \t\tAC\ is also complete.
}

\tetel{\t\tFB.}{Let $X$ be a \bspc\ with the \bap, $\Oa\subset X$ \pscx\ open,
        and $S\to\Oa$ a cohesive sheaf.
         If \pshdom\ is possible in\/ $\Oa$, then\/
\vskip0pt
        {\rm(a)} the cohesive sheaf $S$ admits a complete resolution over\/ $\Oa$,
        and
\vskip0pt
        {\rm(b)} the sheaf cohomology group $H^q(\Oa,S)$ vanishes for all $q\ge1$.
}

\biz{Proof.}{\WLOG\ we may assume by \t\tEA\ that $X$ has a bimonotone \sbs.
	 An inspection of the proof of the analogous \t9.1 in [\rLP] reveals that
        therein it is enough to have \pshdom\ in $\Oa$ and in those subsets of $\Oa$
        to which \p\tEE\ applies, and thus in which \pshdom\ holds by \t\tFA.
         The proof of \t\tFB\ is complete.
}

\biz{Proof of \t\tAD.}{Parts (a) and (b) follow directly from \t\tFB, (c) from
	[\rLP, \S\,10] and \t\tFB, while (d) follows from [\rPB, \t1.3(f)], 
	completing the proof of \t\tAD.
}

\biz{Proof of \t\tAE.}{As the $\delbar$-equation $\delbar u=f$ can be solved locally on balls in $\Oa$ 
	by a \th\ of Defant and Zerhusen [\rDZ]
	(based upon the earlier work [\rLA] of Lempert) a standard step in one
	of the usual proofs of the Dolbeault isomorphism together with \t\tAD(c)
	completes the proof of \t\tAE.
}

	 Further applications of \Th{}s \tAD\ and \tAE\ can also be made, e.g., as
	in [\rDPV] or [\rLP].

\vskip0.20truein
\centerline{\scVIII References}
\vskip0.15truein
\baselineskip=11pt
\parskip=1pt
\frenchspacing
{\rmVIII

	[\rDZ] Defant, A., Zerhusen, A.\,B.,
	{\itVIII
	Local solvability of the $\scriptstyle\delbar$-equation on $\scriptstyle{\Cal L}_1$-spaces},
	Arch. Math.,
	{\bfVIII 90},
	(2008),
	545--553.

	[\rDPV] Dineen, S., Patyi, I., Venkova, M.,
	{\itVIII 
	Inverses depending holomorphically on a parameter in a Banach space},
	J. Funct. Anal., 
	{\bfVIII 237} (2006), no. 1, 338--349.

	[\rLA] Lempert, L.,
	{\itVIII
	The Dolbeault complex in infinite dimensions, II},
	J. Amer. Math. Soc.,
	{\bfVIII 12} (1999), 775--793.

	[\rLB] \vonal,
	{\itVIII
	Approximation of holomorphic functions of infinitely many variables, II},
	Ann. Inst.  Fourier Grenoble
	{\bfVIII 50} (2000), 423--442.

	[\rLC] \vonal,
	{\itVIII
	Plurisubharmonic domination},
	J. Amer. Math. Soc.,
	{\bfVIII 17}
	(2004),
	361--372.

	[\rLP] \vonal, Patyi,~I.,
	{\itVIII
	Analytic sheaves in Banach spaces},
	Ann. Sci. \'Ecole Norm. Sup., 
	4e s\'erie,
	{\bfVIII 40} (2007), 453--486.

	[\rM] Meylan, F.,
	{\itVIII
	Approximation of holomorphic functions in Banach spaces 
	admitting a Schauder decomposition},
	Ann. Scuola Norm. Sup. Pisa, 
	(5) {\bfVIII 5} (2006), no. 1, 13--19.

	[\rPA] Patyi, I.,
	{\itVIII
	On complex Banach submanifolds of a Banach space},
	Contemp. Math., 
	{\bfVIII 435} (2007), 343--354.
			
	[\rPB] \vonal,
	{\itVIII
	On holomorphic Banach vector bundles over Banach spaces},
	Math. Ann.,
	{\bfVIII 341} (2008), no. 2, 455--482.

        [\rP] {Pe{\l}czy{\'n}ski, A.},
        {\itVIII
        Projections in certain {B}anach spaces},
        Studia Math.,
        {\bfVIII 19}
        (1960),
        209--228.
}
\vskip0.20truein
\centerline{\vastag*~***~*}
\vskip0.15truein
{\scVIII
	Imre Patyi,
	Department of Mathematics and Statistics,
	Georgia State University,
	Atlanta, GA 30303-3083, USA,
	{\ttVIII ipatyi\@gsu.edu}}
\bye